\documentclass[a4paper,12pt]{amsart}
\usepackage{amsmath,amsthm,amsfonts,amssymb,amscd,mathrsfs}
\usepackage{psfrag,graphicx}
\usepackage{bbding}
\usepackage{epic,eepic}
\usepackage{color}
\usepackage{ebezier}

\usepackage{epsfig}

\theoremstyle{plain}
\newtheorem{main}{Theorem}

\newtheorem{Thm}{Theorem}[section]
\newtheorem{Lem}[Thm]{Lemma}
\newtheorem{Prop}[Thm]{Proposition}

\theoremstyle{remark}
\newtheorem{Def}[Thm] {Definition}
\newtheorem{Rem}[Thm] {Remark}

\long\def\begcom#1\endcom{}

\newcommand{\length}{\operatorname{\length}}

\newcommand{\dist}{\operatorname{dist}}

\long\def\begcom#1\endcom{}

\def\length{\operatorname{length}}

\def\length{\operatorname{length}}

\begin{document}

\title[Invariant Measures Supported on Chain Recurrent Classes Are Physical-like ]
       {Physical-like Measures Coincide with Invariant Measures Supported on Chain Recurrent Classes}


\author[X. Tian] {Xueting Tian}
\address[X. Tian]{School of Mathematical Science,  Fudan University\\Shanghai 200433, People's Republic of China}
\email{xuetingtian@fudan.edu.cn}
\urladdr{http://homepage.fudan.edu.cn/xuetingtian}

\maketitle

\def\abstractname{\textbf{Abstract}}

\begin{abstract}\addcontentsline{toc}{section}{\bf{English Abstract}} \small
For $C^0$ generic continuous maps or homeomorphisms on compact Riemannian manifold, we prove that (1) the space of physical-like measures coincides with the set of invariant measures supported on chain recurrent classes, (2) every point in the base space is typical (that is, for any point in the base space, its empirical measures  are contained in the space of physical-like measures) and (3) there is a subset of strongly regular set with Lebesgue zero measure but has infinite topological entropy. Moreover, some comparison between $C^0$ generic systems and $C^0$ conservative generic systems are discussed.

\end{abstract}

\section{Introduction}


Let $f : M \rightarrow M $ be a continuous map on a compact   manifold $M$, which does not necessarily preserve any smooth measures with respect to the Lebesgue measure.
Let $\mathcal{P}$ denote the space of all the probability measures endowed with the weak$^*$ topology, and ${\mathcal P}_f \subset {\mathcal P}$  denote the space of $f$-invariant probability measures.
 \label{DefinitionEmpiricalProba}
For any point $x\in M$ and for any integer number $n \geq 1$, the \em empirical probability or time-average measure \em $\Upsilon_n(x)$  of the $f$-orbit of $x$ up to time $n$, is defined by
$$\Upsilon_n(x):= \frac1n\sum_{j=0}^{n-1}\delta_{f^j(x)},$$
where $\delta_y$ is the Dirac  probability measure supported at $y\in M$. Consider the sequence  $\big\{ \Upsilon_n \big\}_{n \in \mathbb{N}^+} $ of empiric probabilities in the space ${\mathcal P}$, and define the  \em  p-omega-limit set \em  $p \omega_f(x) \subset {\mathcal P}$ as follows:
$$p \omega_f(x) := \Big\{\mu \in {\mathcal P}: \ \ \exists \ n_i \rightarrow + \infty \mbox{ such that }
 {\lim_{i \rightarrow + \infty}}\ \Upsilon_{n_i} = \mu \Big\}.$$
It is standard to check that $p \omega_f(x) \subset {\mathcal P}_f$. From \cite{DGS} we know that $p \omega_f(x)$ is always nonempty, weak$^*$-compact and connected.
Among the most useful concepts in the ergodic theory,  the physical probability measures play  an important role.  

\begin{Def} \label{DefinitionPhysical} 
 We call a measure $\mu\in \mathcal P$ {\it physical or SRB} (Sinai-Ruelle-Bowen), if  the set \begin{equation}\label{BasinOfMu}G_\mu=\{x\in M\colon \ p\omega_f(x)=\{\mu\}\}\end{equation}
 has positive Lebesgue measure. The
set $G_\mu$ is called \em basin of statistical attraction \em of $\mu$, or in brief, basin of $\mu$ (even if $\mu$ is not physical).
\end{Def}
\begin{Rem} The above definition of physical or SRB measures is not adopted by all the authors. Some
mathematicians  require  the measure $\mu$ to be ergodic  to call it physical.
Besides, some mathematicians when studying $C^{1 + \alpha}$ systems do not define SRB as a synonym of physical measure, but take into account the property of absolute continuity on the unstable foliation. But, in the scenario of continuous systems, and even for $C^1$ systems, the unstable conditional measures can not be defined because the unstable foliation may not exist (see \cite{Pugh}).

\end{Rem}

Not any system, principally in the $C^1$ context,   possesses physical measures.   This problem can be easily dodged by substituting the definition of physical measure by a weaker concept: \em physical-like \em measure, also called SRB-like or observable measure (see Definition \ref{definitionSRB-like}). 
 In   \cite{CE} it  was proved that for any $C^0$  system $f$,  there exists a nonempty set   composed by all the observable or physical-like measures. 

 \begin{Def}  
  \label{definitionSRB-like}
   Choose any metric $\dist^*$ that induces the weak$^*$ topology on the space ${\mathcal P}$ of probability measures.
A probability measure $\mu\in \mathcal{P}$ is called \em physical-like \em   (or  \em SRB-like \em or \em observable\em) if for any $\varepsilon>0$ the set \begin{equation}\label{BasinEpsilonOfMu} G_\mu(\varepsilon)=\{x\in M\colon \ \mbox{dist}^*(p\omega_f(x),\mu)<\varepsilon\}\end{equation} has positive Lebesgue measure. The set $G_\mu(\varepsilon)$ is called \em basin of $\varepsilon$-partial statistical attraction   \em of $\mu$, or in brief, $\varepsilon$-basin of $\mu$.

\end{Def}


 We denote by $\mathcal{O}_{f}$ the set of  physical-like measures for $f$. It is standard to check that every physical-like measure is $f$-invariant and that ${\mathcal O}_f$ does not depend on the choice of the metric in ${\mathcal P}$. There is a basic fact for ${\mathcal O}_f$.

\begin{Thm}\label{SRB-like} \label{theoremCE}  {\bf (Characterization of physical-like measures \cite{CE})}

Let $f : M \rightarrow M $ be a continuous map on a compact   manifold $M$.
 Then, the set $\mathcal{O}_{f}$ of  physical-like measures 
  is nonempty, weak$^*$ compact, and  contains the limits of the convergent subsequences of the empiric probabilities for Lebesgue almost all the initial states $x \in M$. Besides, no proper subset of ${\mathcal O}_f$ has the latter three properties simultaneously.
\end{Thm}

\vspace{0cm}




A point $x \in M$  is \em irregular \em (see  \cite{Pesin,PP,Barreira-Schmeling2000,Bar,Bar2011}) if the sequence of time-averages along its orbit is not convergent, that is, $\# pw_f(x)>1$. It is also called points  with historic behaviour \cite{Ruelle01,Takens}. Otherwise, $x$ is called \em regular \em (called quasi-regular   in \cite{Oxt,DGS}). Let $IR_f$ be the set of irregular points and $R_f$ be the set of regular points.
In \cite{Barreira-Schmeling2000} it is proved that $IR_f$ carries full topological entropy for hyperbolic systems. This result is generalized to systems with specification-like properties \cite{To2010,Tho2012}. 
The set of strongly regular points (called regular   in \cite{Oxt}), denoted by $SR_f,$ which means that
    $$SR_f=\bigcup_{\mu\in\mathcal P_f,\,\mu\,ergodic} (G_\mu\cap S_\mu),$$ where $S_\mu$ denotes the support of $\mu.$ This set has full measure for any invariant measure by Birkhoff ergodic theorem and ergodic decomposition theorem.




A  $\delta$-pseudo-orbit (or $\delta$-chain) of $f$ from $x$ to $y$ is a sequence $\{x_i\}^n_{i=0}$
with $x_0 = x, x_n = y$ and
$d(f(x_k), x_{k+1}) <\delta$   for $ k = 0, 1, ..., n - 1.$
A point $x \in X$ is called chain recurrent if there is an $\delta$-chain from $x$ to itself for any positive $\delta$. We
can define an equivalence relation on the set of chain recurrent points in such a way that two points x and y are said to be equivalent if for every  $\delta>0$ there
exist an $\delta$-chain from $x$ to $y$ and an $\delta$-chain from $y$ to $x$. The equivalence classes
of this relation are called chain recurrent classes (or chain components).
 These are compact invariant sets
and cannot be decomposed into two disjoint compact invariant sets, hence serve
as building blocks of the dynamics. The topology of
chain recurrent classes and the set of chain recurrent points have been always in particular interest \cite{Ath,RW,Sak,WGW}.




A point $x\in M$ is called {\it physically-typical,} if $ pw_f(x)\subseteq \mathcal O_f.$ Let $h_{top}(A)$ denote the topological entropy of $A\subseteq M$ defined by Bowen \cite{Bowen} and let $Leb(A)$ denotes the Lebesgue measure of $A$. Let $\mathcal{H}(M), \ \mathcal C(M)$ denote the space of all homeomorphisms and continuous maps respectively.  Now we are ready to state our main result.

\begin{main}
 \label{TheoremMain2018}  There is a dense $G_\delta$ subset $\mathcal{R}\subseteq \mathcal{H}(M)\ (or \ \mathcal C(M))$ in $C^0$ topology such that for any  $f\in \mathcal R,$    \\
 {
 (1)    $ \mathcal O_f
 =\{\mu|\ S_\mu   \text{ is contained in  a chain recurrent class}\}.$
   \\
 (2) for any $x\in M$, $x$ is physically-typical.  \\ 
 (3) there is a subset $\Lambda\subseteq SR_f$ such that $$Leb(\Lambda)=0\text{ but }h_{top}(\Lambda)=h_{top}(f)=+\infty.
$$}
 \end{main}

Here we also list a  recent related result \cite{DOT} that there is a dense $G_\delta$ subset $\mathcal{R}\subseteq \mathcal{H}(M)\ (or \ \mathcal C(M))$ such that for any  $f\in \mathcal R,$    \\ 
(4) there is a subset $\Lambda\subseteq R_f$ such that $$Leb(\Lambda)=1\text{ but }h_{top}(\Lambda)=0;\  \ 
$$
(5)  $$Leb(IR_f)=0 \text{ but }h_{top}(IR_f)=h_{top}(f)=+\infty.
$$ Here item (3) of Theorem \ref{TheoremMain2018} gives us more information on Lebesgue measure and topological entropy in  $SR_f$ which still  has full measure for any invariant measure. 

\begin{Rem} By \cite{CTrou}  for a $C^0$ generic $f,$ periodic measures are dense in $\mathcal O_f$ so that $\inf_{\mu \in \mathcal O_f}h_\mu(f)=0.$ For $C^1$ partially hyperbolic diffeomorphisms (which form an open set in the space of diffemorphisms), it is proved in \cite{CT2016} that $\inf_{\mu \in \mathcal O_f}h_\mu(f)>0,$ where $h_\mu(f)$ denotes the metric entropy of $\mu.$ Thus $\{f:\inf_{\mu \in \mathcal O_f}h_\mu(f)=0\}$ is not $C^r$ generic for $r\geq 1.$ For $C^1$ partially hyperbolic diffeomorphisms, note that $\inf_{\mu \in \mathcal O_f}h_\mu(f)>0$ implies that any periodic measure (if exists) is not in $\mathcal O_f$  and any periodic point (if exists) is not physically-typical. Combining with classical Kupka-Smale's theorem that periodic measures and points exist and must be hyperbolic in $C^1$ generic diffeomorphisms,  items (1) and (2) of Theorem \ref{TheoremMain2018} both are false in $C^1$ generic diffeomorphisms.

\end{Rem}

\begin{Rem}

A point $x\in M$ is called without physical-like behavior if  $pw_f(x)\cap \mathcal O_f=\emptyset.$ Denote $\Gamma_f$ be the set of such points. By Theorem \ref{theoremCE}, $\Gamma$ always has zero Lebesgue measure.   It is proved in \cite{CTV} that such points form a nonempty set with full topological entropy in all transitive Anosov diffeomorphisms for which such systems form an open subset in the space of diffeomorphisms by structural stability. This implies that $\{f|\ \Gamma_f= \emptyset\}$ is not  $C^r$ generic for any $r\geq 1$,  different from item (2) of Theorem \ref{TheoremMain2018} which implies that $\{f|\ \Gamma_f= \emptyset\}$ is $C^0$ generic. 
\end{Rem}

\begin{Rem} Takens' last problem said that: Whether there are persistent classes of smooth
dynamical systems such that the set of initial states which give rise to orbits with historic
behaviour has positive Lebesgue measure? S. Kiriki,  T. Soma  \cite{KS} proved that any Newhouse open set in the space of $C^r$-diffeomorphisms ($r\geq 2$) on a closed surface is contained in the closure of the set of diffeomorphisms which have non-trivial wandering domains whose forward orbits have historic behavior so that Takens' last problem holds. This implies that  the parts for Lebesgue measure in item (4) and  item (5) are false for $C^r$ generic diffeomorphisms ($r\geq 2$).

\end{Rem}

We point out some similar results as Theorem \ref{TheoremMain2018} and \cite{DOT} for $C^0$ generic conservative systems. Let $\mathcal{H}_{Leb}(M), \ \mathcal C_{Leb}(M)$ denote the space of all conservative  homeomorphisms and continuous maps respectively.
\begin{Thm} \label{Thm-1-same} There is a dense $G_\delta$ subset $\mathcal{R}\subseteq \mathcal{H}_{Leb}(M)\ (or \ \mathcal C_{Leb}(M))$ in $C^0$ topology such that for any  $f\in \mathcal R,$  above item (3), (4) and (5) are true.

\end{Thm}

 On the other hand,  item (2) of Theorem \ref{TheoremMain2018} implies that for $C^0$ generic   systems, $\Gamma_f$ is  empty  (carrying zero topological entropy) and $M\setminus \Gamma_f=M$ (carrying full topological entropy), but here we point out  a  following different result  for $C^0$ generic conservative systems.

\begin{Thm} \label{Thm-2-different} There is a dense $G_\delta$ subset $\mathcal{R}\subseteq \mathcal{H}_{Leb}(M)\ (or \ \mathcal C_{Leb}(M))$ such that for any  $f\in \mathcal R,$  above item (1) and (2) of Theorem \ref{TheoremMain2018} are false; and moreover
\\
(6)  $\Gamma_f$ is not only nonempty but also carries full infinite  topological entropy and $M\setminus \Gamma_f\subsetneq M$ carries zero topological entropy;\\
(7) $Leb(IR_f)=0$ and $h_{top}(IR\cap \Gamma_f)=h_{top}(f)=+\infty;$\\ 
 (8) 
  there is a subset $\Lambda\subseteq SR_f$ such that $$Leb(\Lambda)=0\text{ but }h_{top}(\Lambda\cap \Gamma_f)=h_{top}(f)=+\infty;$$
 (9)  there is a subset $\Lambda\subseteq R_f\setminus SR_f$ such that $$Leb(\Lambda)=0\text{ but }h_{top}(\Lambda\cap \Gamma_f)=h_{top}(f)=+\infty.$$

\end{Thm}
\begin{Rem} Note that items (6), (7), (8) and (9) are false for $C^0$ generic systems, since by item (2) of Theorem \ref{TheoremMain2018} one has $\Gamma_f=\emptyset.$  Note that item (3) is weaker than item (8) and  item (5) is weaker than item (7). Similarly one can ask a weaker version of item (9) for a $C^0$ generic system, whether there is a subset $\Lambda\subseteq R_f\setminus SR_f$ such that $$Leb(\Lambda)=0\text{ but }h_{top}(\Lambda)=h_{top}(f)=+\infty?$$ Up to now this is still unknown.

\end{Rem}

\begin{Rem} By  \cite{CTrou} for $C^0$ generic systems,  ergodic measures are dense in $\mathcal O_f$ so that by variational principle  $\sup_{\mu \in \mathcal O_f}h_\mu(f)=h_{top}(f)=+\infty.$  However, for $C^0$ generic conservative systems, $$\sup_{\mu \in \mathcal O_f}h_\mu(f)=0,$$ since in this case Lebesgue measure is ergodic with zero metric entropy \cite{Gui} and by Theorem \ref{theoremCE} $\mathcal O_f$ is  just composed by the Lebesgue measure.

\end{Rem}

The proofs of these results are not only based on some classical results but also based on some very recent results.

\section{Periodic Approximation} {  Let $(X,f)$ be a topological dynamical system, meaning that $f:X\rightarrow X$ is a continuous map on a compact metric space $X$.}
\subsection{Various periodic orbit tracing properties}
In this section we recall some orbit tracing properties. For more details and results, one may refer to ~\cite{DGS,
Sig
}.
  For any $\delta>0$,   a sequence $\{x_n\}_{n=0}^{+\infty}$ is called a \textit{$\delta$-pseudo-orbit} if
$$d(f(x_n),  x_{n+1})<\delta~~\text{for any}~~n\geq 1.  $$
A $\delta$-pseudo-orbit $\{x_n\}_{n=0}^{+\infty}$ is called a {\it periodic $\delta$-pseudo-orbit}, if there is $T\in\mathbb{N}$ such that $x_{i+T}=x_i, \forall i\geq 0,$ where the smallest positive integer $T$ satisfying this is called its {\it period}.
A pseudo-orbit $\{x_n\}_{n=0}^{+\infty}$ is \textit{$\epsilon$-shadowed} by the orbit of a point $y\in X$, if
$$d(f^n(y),  x_n)<\epsilon \ ~\textrm{for any}~~n\in\geq 1.  $$

\begin{Def}\label{Def:shadowing}
We say that $(X, f)$ satisfies the \textit{(periodic) shadowing property}, if for any $\epsilon>0$,   there exists  $\delta>0$ such that any (periodic) $\delta$-pseudo-orbit is $\epsilon$-shadowed by the orbit of a (periodic) point in $X$.
\end{Def}


We introduce a notion of gluing property on a set for which the tracing point may be not in this set.



\begin{Def}\label{Def:gluing} Let $\Lambda\subseteq X$ be an invariant and closed subset.
We say that $\Lambda$ has the \emph{periodic gluing orbit property}, if for any $\varepsilon>0$ there exists $M(\varepsilon)\in\mathbb{N}$ such that
for any points $x_1,\cdots ,x_k\in X$ and any integers $n_1,\cdots ,n_k\geq 1,$ there exist $m_1,\cdots ,m_{k}\leq M(\varepsilon)$ and a periodic point $x\in X$  with period $\sum_{j=1}^{k} (m_j+n_j)$  satisfying
$d(f^{l+\sum_{j=0}^{i-1} (m_j+n_j)}(x),f^l(x_i))<\varepsilon$ for all $0\leq l\leq n_i-1$ and $1\leq i \leq k,$ where $m_0=n_0=0.$
\end{Def}


\begin{Lem}\label{Lem-gluing} Suppose $(X,f)$ has periodic shadowing property and $\Lambda$ is a chain recurrent class. Then $\Lambda$ has periodic gluing property.

\end{Lem}

{\bf Proof. } 
For any
$\varepsilon>0, $ 
 by periodic shadowing there
exists $\delta>0$ such that any periodic $\delta$-pseudo orbit in
$\Lambda$ can be $\varepsilon$ shadowed by a periodic orbit in $M$.

Take and fix for $\Lambda$ a finite cover
$\alpha=\{U_1,U_2,\cdots,U_{r_0}\}$ by nonempty open balls $U_i$ in
$\Lambda$ satisfying $diam(U_i) <\delta$, $i=1,2,\cdots,r_0$. Since
$\Lambda$ is chain transitive, for $x\in U_i, y\in U_j,$ and any $i,j=1,2,\cdots,r_0,$ there exist a
positive integer $X_{i,j}$ and $\delta$-chain $\Gamma_{i,j}:=\{x^0_{i,j},x^1_{i,j},\cdots, x^{X_{i,j}}_{i,j}\}$
such that $$x^0_{i,j}=x\in U_i,\, x^{X_{i,j}}_{i,j}=y\in U_j. $$
 Let
$$M(\epsilon):=max_{1\leq
i\neq j\leq r_0}X_{i,j}.$$

Now let us consider a given  sequence of points
$x_1,x_2,\cdots,x_N\in\,\Lambda,$ and a sequence of positive numbers
${n_1,n_2,\cdots,n_N}$. Take and fix $U_{i_0},U_{i_1}\in\alpha$ so
that $x_i\in U_{i_0},f^{n_i}x_i\in U_{i_1},i=1,2,\cdots,N.$ Take  $\delta$-chain $\Gamma_{i_1,{(i+1)_0}}:=\{x^0_{i_1,(i+1)_0},x^1_{i_1,(i+1)_0},\cdots, x^{X_{i_1,(i+1)_0}}_{i_1,(i+1)_0}\}$   such that $$x^0_{i_1,(i+1)_0}=f^{n_i}x_i\in U_{i_1},\, x^{X_{i_1,(i+1)_0}}_{i_1,(i+1)_0}=x_{i+1}\in U_{(i+1)_0}  $$ for
$i=1,2,\cdots,N$ where $x_{N+1}=x_1,U_{(N+1)_0}=U_{1_0}$   Thus we get a periodic
 $\delta$-pseudoorbit:
$$\{f^t(x_1)\}_{t=0}^{n_1}\cup
\Gamma_{1_1,2_0}\cup \{f^t(x_2)\}_{t=0}^{n_2}
\cup\cdots\cup\{f^t(x_N)\}_{t=0}^{n_N}
\cup\Gamma_{N_1,(N+1)_0}.$$  Hence there
exists a periodic point $z\in M$ $\varepsilon$-shadowing the above
sequence with period $p=\sum_{i=1}^{N}[n_i+X_{i_1,(i+1)_0}].$
Clearly
$p\in[\sum_{i=1}^{N}n_i,\,\,\sum_{i=1}^{N}n_i+NM(\epsilon)].$
More precisely,
$$d(f^{c_{i-1}+j}z,f^jx_i)<\varepsilon,\,\,j=0,1,\cdots,n_i,\,\,\,i=1,2,\cdots,N,$$
where $c_i$ is defined
 as follows:
 $$c_i=\begin{cases}
 \,\,0,&\text{for }\,\,i=0\\
 \,\,\sum_{j=1}^{i}[n_j+X_{j_1,(j+1)_0}],\,\,&\text{for}\,\,i=1,2,\cdots,N.\\
\end{cases}
 $$ \qed 

\subsection{Approximation by periodic measures}

Let $\mathcal P_p(f)$ denotes the space of all periodic measures.
\begin{Lem}\label{Lem-densityPeriodic} Suppose   $\Lambda$ is a compact invariant set with periodic gluing property. Then $\mathcal P_f(\Lambda)\subseteq \overline{\mathcal P_p(f)}.$

\end{Lem}
{\bf Proof. } Take a sequence of non-zero functions $\{f_n\}$ dense in the space $C^0(X,\mathbb{R})$ of all continuous functions. Then the weak$^*$ topology can be defined as following: for any $\mu,\nu\in \mathcal P(X),$
$$dist(\mu,\nu):=\sum_{n=1}^\infty\frac{|\int g_n d\mu-\int g_n d\nu|}{2^{n+1}\|g_n\|}.$$
 Fix $\epsilon >0$ and $L>0$ such that $$\sum_{n=L+1}^\infty\frac{2}{2^{n+1} }<\frac{\epsilon}{5}.$$
 Let $ \mu\in \mathcal P_f(\Lambda)$ and $F=\{\frac{g_i}{\|g_i\|}\}_{i=1}^L$. Then one has to show that there is $\mu_p\in \mathcal P_p(X)$ such that for any $g\in F,$ $|\int g d\mu-\int g d\mu_p|<\frac{4\epsilon}{5}. $

Choose $\delta>0$ such that for all $g\in F$ one has $|g(x)-g(y)|<\frac{\epsilon}{5}$ whenever $d(x,y)<\delta,\ x,y\in X.$ Let $M(\delta)$ be the number given in the gluing property.

 By ergodic decomposition theorem for the subsystem $f|_\Lambda$, there exist finite ergodic measures $\mu_i\in \mathcal P(\Lambda)$ ($i=1,2\cdots,k$) and positive numbers $\theta_i$ with $\sum_{i=1}^{k}\theta=1$  such that for any $g\in F,$
 $$|\int g d\mu-\sum_{i=1}^{k}\theta_i\int g d\mu_i |=|\int g d\mu-\int g d(\sum_{i=1}^{k}\theta_i\mu_i) |<\frac{\epsilon}{5}. $$
 We may assume $\theta_i$ are all rational numbers. Then there exist natural numbers $m_i$ such that $\theta_i=\frac{m_i}{m}$ where $m=\sum_{i=1}^{k}m_i.$
 By Birkhoff ergodic theorem for the subsystem $f|_\Lambda$, we can choose $x_i\in \Lambda$ such that $$\lim_{n\rightarrow\infty}\frac{1}{n} \sum_{j=1}^{n-1} g(f^jx_i)=\int g d\mu_i,\,\,i=1,2,\cdots,k,\ g\in F.$$
Take $N$ large enough such that $\frac{2M(\delta)}{N}<\frac{\epsilon}{5}$ and $$ |\frac{1}{N} \sum_{j=1}^{N-1} g(f^jx_i)-\int g d\mu_i| <\frac{\epsilon}{5},\,\,i=1,2,\cdots,k,\ g\in F.$$

Let $\{y_n|1\leq n\leq mN\}\subseteq \Lambda$  be  a sequence of points defined by letting $y_n$ run $m_1$ times through $\{x_1, fx_1,\cdots,f^{N-1}x_1\},$ then $m_2$ times through $\{x_2, fx_2,\cdots,f^{N-1}x_2\}, $
etc., finally $m_k$ times through $\{x_k, fx_k,\cdots,f^{N-1}x_k\}.$
 Then for any $g\in F,$
 $$|\int g d\mu-\frac1{mN}\sum_{i=1}^{mN}  g (y_i) |<\frac{2\epsilon}{5}. $$ By gluing property, there exist    $0\leq q_i\leq M(\delta), i=1,2,\cdots,m,$   and a periodic point $z\in X$  with period $p=mN+\sum_{j=1}^{m} q_j$  satisfying
$d(f^{t}(z),f^t(y_i))<\varepsilon$ for all $a_t\leq t\leq b_t, \ t=0,1\cdots,m-1$   where $a_t=tN+\sum_{i=1}^{t}q_i, \ b_t=tN+\sum_{i=1}^{t}q_i+N.$

 Note that for any $g\in F,$ $$|\frac1p \sum_{i=1}^{p}g (f^iz)-\frac1{mN}\sum_{i=1}^{mN}  g (y_i) |\leq \frac{2M(\delta)}{N}+\frac{\epsilon}{5}<\frac{2\epsilon}{5}. $$ Thus for any $g\in F,$  $$|\int g d\mu-\frac1p \sum_{i=1}^{p}g (f^iz)|<\frac{4\epsilon}{5}. $$ Now we complete the proof.
\qed



\subsection{Basic facts}

Let $A\subseteq X$ be a nonempty compact invariant set.   We call $A$ \emph{internally chain transitive} if for any $a,  b\in A$ and any $\epsilon>0$,   there is an $\epsilon$-chain  $\{x_i\}^n_{i=0}$  contained in $A$ with connecting $a$ and $b$.
Let $\omega_f=\{A\subseteq X:\exists x\in X~\text{with}~\omega_f(x)=A\}$ and denote the collection of internally chain transitive sets by $ICT.$
It has been shown  \cite{HSZ} that every $\omega$-limit set is internally chain transitive,
and the converse has also been shown in a variety of contexts, including Axiom A
diffeomorphisms \cite{Bowen1975}, shifts of finite type  \cite{BGKR}, topologically hyperbolic maps \cite{BGOR}, and in
certain Julia sets  \cite{BR,BMR}, systems with orbit limit shadowing \cite{GM}.
  From \cite{MR} for systems with shadowing property, $ \overline{\omega_f}= ICT.$

{
\begin{Prop}\label{Prop-physicallike-chainrec} For any dynamical system $(X,f), $ \\
(1) $\{\mu|\ S_\mu \subseteq \Lambda \text{ where }\Lambda \in  {\omega_f}\}\subseteq   \{\mu|\ S_\mu \subseteq \Lambda \text{ where }\Lambda \in \overline{\omega_f}\}\subseteq \{\mu|\ S_\mu \subseteq \Lambda \text{ where }\Lambda \in ICT\}\subseteq \{\mu|\ S_\mu \subseteq \Lambda \text{ where }\Lambda \text{ is a chain recurrent class}\}.$ \\
(2) for any $x\in M,$ $ pw_f(x)\subseteq   \{\mu|\ S_\mu \subseteq \Lambda \text{ where }\Lambda \in  {\omega_f}\}
.$
\\
(3) If $X$ is a compact manifold, then $$\mathcal O_f\subseteq \{\mu|\ S_\mu \subseteq \Lambda \text{ where }\Lambda \in \overline{\omega_f}\} \cap \overline{\{\mu|\, \mu\in pw_f(x), x\in M\}}. 
$$
(4) If $X$ is a compact manifold and $f$ has periodic shadowing, then $\mathcal O_f\subseteq
\overline{\mathcal P_p(f)}.$ \\
\end{Prop}

{\bf Proof.}  (1) From \cite{MR}  we know that   $ICT$ is closed in Hausdorff distance. By \cite{HSZ} $  {\omega_f}\subseteq ICT$ and thus $ \overline{\omega_f}\subseteq ICT.$ Note that any internally chain transitive set should be contained in one chain recurrent class. Thus one has item (1).

(2) It is easy to check for any $x\in X$, $\cup_{\mu\in pw_f(x)}S_\mu\subseteq \omega_f(x)$ and thus for any $\mu\in pw_f(x)$, $\mu(\omega_f(x))=1.$ Thus one has item (2).

(3)
 Given $\mu\in \mathcal O_f, $ by definition there exists
$x_n\in M$ such that $dist^*(pw_f(x_n),\mu)<\frac1n.$  By the choice of $x_n,$ there exists $\mu_n\in pw_f(x_n)$ such that $d(\mu_n,\mu)<\frac1n$ and $\mu_n(\omega_f(x_n))=1.$
Let $\Lambda_n=\omega_f(x_n)$. Take a convergent subsequence $\{ \Lambda_{n_k}\}$ in Hausdorff distance for which the limit set is a closed set, denoted by $\Lambda.$ Then $\Lambda\in \overline{\omega_f}.$ For $\epsilon>0,$ let $\Lambda_\epsilon:=\{x|\ d(x,\Lambda)\leq \epsilon\}$. Then $\mu(\Lambda_\epsilon)\geq \lim_{n\rightarrow\infty}\mu_{n_k}(\Lambda_\epsilon)=1$ and thus $\mu(\Lambda)=1.$
 Now one ends the proof of item (3).

 (4)  By Proposition \ref{Prop-physicallike-chainrec} (1), (3), Lemma \ref{Lem-gluing} and Lemma \ref{Lem-densityPeriodic}, $$\mathcal O_f\subseteq \{\mu|\ S_\mu   \text{ is contained in  a chain recurrent class}\}\subseteq
\overline{\mathcal P_p(f)}.$$
\qed
}

\section{Proofs}

 \bigskip

\subsection{ Proof of Theorem \ref{TheoremMain2018}.}

(1)  It is known (see \cite{KMOP} and \cite{Kos2005,Kos}, respectively) that periodic shadowing is a
generic property in $\mathcal C(M)$ and in $\mathcal H(M)$.  By Proposition \ref{Prop-physicallike-chainrec} (1), (3), Lemma \ref{Lem-gluing} and Lemma \ref{Lem-densityPeriodic}, $$\mathcal O_f\subseteq \{\mu|\ S_\mu   \text{ is contained in  a chain recurrent class}\}\subseteq
\overline{\mathcal P_p(f)}.$$ By \cite{CTrou} $\mathcal O_f=\overline{\mathcal P_p(f)}$ is a generic property in $\mathcal C(M)$ and in $\mathcal H(M)$.   Thus one ends the proof of item (1) of Theorem \ref{TheoremMain2018}.

(2) By Proposition \ref{Prop-physicallike-chainrec} (1), (2) and Theorem \ref{TheoremMain2018} (1),  for any $x\in M,$ $pw_f(x)\subseteq \mathcal O_f.$ 

(3) Fix $f\in\mathcal{H}(M)\ (or \ \mathcal C(M))$. By classical variational principle, take $\mu_n$ to be a sequence of ergodic measures with $\lim_{n\rightarrow\infty }h_{\mu_n}(f)=h_{top}(f).$
Let $\Lambda=\cup_{i=1}^\infty G_{\mu_i}\cap S_{\mu_i}$. Then $\Lambda\subseteq SR_f.$ Note that $\mu_i( G_{\mu_i}\cap S_{\mu_i})=1$ and thus $\mu_i(\Lambda)=1.$
From \cite[Theorem 1]{Bowen} one knows that for any invariant measure $\omega$ and any set $Y\subseteq M$ with $\omega(Y)=1$, then $h_\omega(f)\leq h_{top}(Y).$ Thus $h_{top}(\Lambda)\geq \lim_{n\rightarrow\infty }h_{\mu_n}(f)=h_{top}(f).$

From \cite{CTrou,AA} there is a dense $G_\delta$ subset $\mathcal{R}_1\subseteq \mathcal{H}(M)\ (or \ \mathcal C(M))$ such that for any  $f\in \mathcal R_1,$ there is no physical measure. Thus $Leb(G_{\mu_i})=0$ so that $Leb(\Lambda)=0.$ From \cite{Yo,Kos} there is a dense $G_\delta$ subset $\mathcal{R}_2\subseteq \mathcal{H}(M)\ (or \ \mathcal C(M))$ such that for any  $f\in \mathcal R_2,$ $h_{top}(f)=+\infty.$  Let $\mathcal{R}=\mathcal{R}_1\cap \mathcal{R}_2$ and then one ends the proof. \qed

\subsection{Proof of Theorem \ref{Thm-2-different}}

Recall a general fact by Bowen \cite{Bowen} that letting $QR(t):=\{x\in X| \text{ there is }\mu \in V_f(x)\text{ such that } h_\mu\leq t\}$, then $h_{top}(T, QR(t))\leq t.$ From \cite{Bowen} we also know that
$$
h_{top}(\cup_{i=1}^\infty Z_i)=\sup_{i=1}^\infty h_{top}(Z_i).$$
 Then it is easy to check that $ M\setminus \Gamma_f\subseteq QR(t)$ when $t=sup_{\mu\in \mathcal O_f}h_\mu(f)$). Thus

\begin{Lem}\label{Lem-O-f-entropy}
For any $A\subseteq M, $
$h_{top}(A\setminus \Gamma_f)\leq h_{top}(M\setminus \Gamma_f)\leq  sup_{\mu\in \mathcal O_f}h_\mu(f).$ If further $sup_{\mu\in \mathcal O_f}h_\mu(f)<h_{top}(f)$ and $h_{top}(A)=h_{top}(f)$, then $h_{top}( A\cap \Gamma_f)=h_{top}(f)$.
\end{Lem}
From \cite{Gui} we know that for a $C^0$ generic conservative system, $f$ is topological mixing with infinite topological entropy, the periodic points are dense in the whole space,  the Lebesgue measure $Leb$ is ergodic with zero metric entropy and thus $\mathcal O_f$ is composed by just the Lebesgue measure. Recently it is proved in \cite{GL} for a $C^0$ generic conservative system, $f$ has specification property (and also has shadowing and periodic shadowing). Note that any periodic measure   is not in $\mathcal O_f$  and any periodic point is not physically-typical.  Thus  items (1) and (2) of Theorem \ref{TheoremMain2018} both are false.

\medskip
(6) Since $sup_{\mu\in \mathcal O_f}h_\mu(f)=h_{Leb}(f)=0,\ h_{top}(f)=+\infty,$ by Lemma \ref{Lem-O-f-entropy} one gets that  $h_{top}(M\setminus \Gamma_f)=0 $ and $h_{top}(   \Gamma_f)=h_{top}(f)=+\infty$.

\medskip
(7) Note that $Leb(G_{Leb})=1$ and $G_{Leb}\subseteq R_f$ so that $Leb(IR_f)=0.$ Given a continuous function $\phi:X\rightarrow \mathbb{R}$, let
$$I_\phi:=\{x\in M|\ lim_{n\rightarrow \infty}\frac 1n \sum_{i=0}^{n-1}\phi(f^ix) \text{ does not converge}\}.$$
 Since $f$ is not uniquely ergodic, then there exists a continuous function $\phi:X\rightarrow \mathbb{R}$ such that $\inf_{\mu\in \mathcal P_f}\int \phi d\mu  <\sup_{\mu\in \mathcal P_f}\int \phi d\mu. $ 
 Then by \cite{Tho2012} one can get that $$h_{top}(I_\phi) =h_{top}(f).$$
Note that $I_\phi\subseteq IR_f$ so that $h_{top}(IR_f) =h_{top}(f).$
Then by 
 Lemma \ref{Lem-O-f-entropy}, one has item (7).
 \medskip

(8) 
Fix small $\epsilon>0.$ By variational principle, there is an ergodic $\mu$ such that $h_\mu(f)>h_{top}(f)-\epsilon.$  By Birkhoff ergodic theorem $\mu(G_\mu \cap S_\mu)=1$ and then by \cite{Bowen} $h_{top}(G_\mu\cap S_\mu)\geq h_\mu(f)>h_{top}(f)-\epsilon.$ Note that $\mu$ is not the Lebesgue measure since $h_{\mu}(f)>0$. Then $G_\mu \cap S_\mu\subseteq SR_f\cap \Gamma_f$ and thus one has item (8).

\medskip
(9) By \cite{PS} specification implies that for any invariant measure $\mu$, $h_{top}(G_\mu)= h_\mu(f). $ Fix small $\epsilon>0.$ By variational principle, there is an ergodic $\mu$ such that $h_\mu(f)>h_{top}(f)-\epsilon.$ Since $f$ is not uniquely ergodic, one can take an invariant measure $\mu'\neq \mu$. Choose $\theta\in(0,1)$ close 1 enough such that $h_{\nu}(f)>h_{top}(f)-\epsilon,$ where $\nu=\theta\mu+(1-\theta)\mu'.$
Then $\nu$ is not ergodic, not the Lebesgue measure, and  $h_{top}(G_\nu)= h_\nu(f)>h_{top}(f)-\epsilon. $  Note that $G_\mu  \subseteq (R_f\setminus SR_f)\cap \Gamma_f$ and thus one has item (9). \qed



\subsection{Proof of Theorem \ref{Thm-1-same}}
By Theorem \ref{Thm-2-different}, item (8) implies item (3) and item (7) implies (5). For item (4),  recall that from \cite{Gui}   for a $C^0$ generic conservative system,   the Lebesgue measure $Leb$ is ergodic with zero metric entropy. Taking $\Lambda=G_{Leb}$, then by ergodicity of $Leb$ and Birkhoff ergodic theorem $Leb(\Lambda)=1$. And by \cite{Bowen} $h_{top}(\Lambda)\leq h_{Leb}(f)=0.$ \qed

\section*{Acknowlegements}
{    X. Tian is  supported by National Natural Science Foundation of China  (grant no. 11671093).


\begin{thebibliography}{10}



\small


\bibitem{AA}F. Abdenur and M. Andersson, {\it Ergodic theory of generic continuous maps,}
Comm. Math. Phys. 318 (2013), no. 3, 831-855.

 \bibitem{Ath} K. Athanassopoulos, {\it One-dimensional chain recurrent sets of flows in the 2-sphere,} Math.
Z. 223(1996), 643-649.


  \bibitem{Bar} L. Barreira, {\it Dimension and recurrence in hyperbolic dynamics}, Progress in Mathematics, vol.
272, Birkh$\ddot{a}$user, 2008.

\bibitem{Bar2011}
L. Barreira, {\it Thermodynamic formalism and applications to dimension theory}. Springer Science $\&$ Business Media, 2011.

\bibitem{Ba-Sau-TAMS}
Barreira L, Saussol B. {\it Variational principles and mixed multifractal spectra}. Transactions of the American Mathematical Society, 2001, 353(10):3919-3944.
\bibitem{Barreira-Schmeling2000} L. Barreira, J. Schmeling, {\it Sets of non-typical points have full topological entropy and full Hausdorff dimension,}  Israel Journal of Mathematics, 2000, 116(1): 29-70.



\bibitem{BGKR} A. D. Barwell, C. Good, R. Knight and B. E. Raines. {\it  A characterization of $\omega$-limit sets in shift spaces.}
Ergod. Th. \& Dynam. Sys. 30(1) (2010), 21-31.


\bibitem{BGOR} A. D. Barwell, C. Good, P. Oprocha and B. E. Raines. {\it Characterizations of  $\omega$-limit sets of topologically
hyperbolic spaces.} Discrete Contin. Dyn. Syst. 33(5) (2013), 1819-1833.


\bibitem{BMR}A. D. Barwell, J. Meddaugh and B. E. Raines. Shadowing and $\omega$-limit sets of circular Julia sets. Ergod. Th.
\& Dynam. Sys. 35(4) (2015), 1045-1055.


\bibitem{BR} A. D. Barwell and B. E. Raines. The $\omega$-limit sets of quadratic Julia sets. Ergod. Th. \& Dynam. Sys. 35(2)
(2015), 337-358.





\bibitem{Bowen71-trans} R. Bowen,  {\it Periodic point and measures for Axiom-A diffeomorphisms,} Trans. Amer. Math. Soc. 154  (1971), 377-397.


\bibitem{Bowen1975} R. Bowen.  {\it $\omega$-limit sets for axiom A diffeomorphisms.} J. Differential Equations 18(2) (1975), 333-339.

\bibitem{Bow} R. Bowen, {\it Periodic orbits for hyperbolic flows,}  Amer. J. Math., 94 (1972), 1-30.

 \bibitem{Bowen2} R. Bowen, {\it Equilibrium states and the ergodic theory of Anosov
diffeomorphisms}, Springer, Lecture Notes in Math. 470  (1975).

\bibitem{Bowen} R. Bowen, {\it Topological entropy for noncompact sets,}  Trans. Amer. Math. Soc. 184 (1973), 125-136.






\bibitem{CE} E. Catsigeras,    H. Enrich,  {\it SRB-like measures  for $C^0$ dynamics,}  Bull. Pol. Acad. Sci. Math. 59, 2011, 151-164.

\bibitem{CT2016}E Catsigeras, X Tian,{\it Dominated Splitting, Partial Hyperbolicity and Positive Entropy},
Discrete and Continuous Dynamical System - Series A 36 (9), 4739-4759.



\bibitem{CTV} E Catsigeras, X Tian, E Vargas, {\it Topological Entropy on Points without Physical-like Behaviour,}
Mathematische Zeitschrift, to appear.


\bibitem{CTrou} E. Catsigeras,    S. Troubetzkoy,  {\it Invariant measures for typical continuous maps on manifolds},  arXiv:1811.04805.

\bibitem{DOT} Dong Y,   Oprocha P,   Tian X.   \emph{On the irregular points for systems with the shadowing property.  } Ergod.   Th.   Dynam.   Sys. ,   2018, 38 (6), 2108-2131.


\bibitem{DongTian2016} Y. Dong, X. Tian, {\it  Different Statistical Future of Dynamical Orbits over Expanding or Hyperbolic Systems (I): Empty Syndetic Center}, arXiv:1701.01910v2.

\bibitem{DGS}
M. Denker, C. Grillenberger and K. Sigmund,  {\it Ergodic Theory on the
Compact Space,} Lecture Notes in Mathematics {\text{527}}.


\bibitem{GM} C.
Good , J.Meddaugh, {\it  Orbital shadowing, internal chain transitivity and ${\it\omega} $-limit sets.}  Ergodic Theory and Dynamical Systems. 2018, 38(1), 143-154.
\bibitem{Gui} P A Guiheneuf, {\it Proprietes dynamiques generiques des homeomorphismes conservatifs
(French, with English and French summaries),} Ensaios Matematicos [Mathematical
Surveys], vol. 22, Sociedade Brasileira de Matematica, Rio de Janeiro, 2012.

\bibitem{GL} Guiheneuf P A, Lefeuvre T. {\it On the genericity of the shadowing property for conservative homeomorphisms.}  Proceedings of the American Mathematical Society. 2018; 146(10):4225-37.

\bibitem{HSZ} M. W. Hirsch, H. L. Smith and X.-Q. Zhao. {\it  Chain transitivity, attractivity, and strong repellors for
semidynamical systems.} J. Dynam. Differential Equations 13(1) (2001), 107-131.


\bibitem{KS}
S. Kiriki and T. Soma, \emph{Takens' last problem and existence of non-trivial wandering domains},  Advances in Mathematics, 306,  524-588, 2017.


\bibitem{Kos2005} P. Koscielniak. {\it On genericity of shadowing and periodic shadowing property. } J. Math. Anal. Appl. 310
(2005), 188-196.

\bibitem{Kos} P. Koscielniak. {\it On the genericity of chaos,} Topology Appl. 154 (2007), 1951-1955.


\bibitem{KMOP} P. Koscielniak, M. Mazur, P. Oprocha and P. Pilarczyk., {\it Shadowing is generic-a continuous case.} Discrete
Contin. Dyn. Syst. 34 (2014), 3591-3609.





\bibitem{MR} J. Meddaugh and B. E. Raines. {\it Shadowing and internal chain transitivity.} Fund. Math. 222 (2013), 279-287.




{ \bibitem{Oxt}  J. C. Oxtoby, {\it Ergodic sets,}  Bull. Amer. Math. Soc. 58, 116-136 (1952).}



\bibitem{Pesin}
Y. B. Pesin, \emph{Dimension theory in dynamical systems: contemporary views and applications}, Chicago Lectures in Mathematics, University of Chicago Press, 2008.

\bibitem{PP}
Y. B. Pesin and B. S. Pitskel$'$, \emph{Topological pressure and the variational principle for noncompact sets}, Functional Analysis and its Applications, 1984, 18(4): 307-318.


\bibitem{PS2005} C.-E. Pfister, W.G. Sullivan, {\it Large Deviations Estimates for Dynamical Systems
without the Specification Property. Application to the $\beta$-shifts,} Nonlinearity 18,
237-261 (2005).

\bibitem{PS} C. Pfister, W.  Sullivan,
{\it  On the topological entropy of saturated sets,} Ergod. Th. Dynam. Sys.
27, 929-956 (2007).


\bibitem{Pugh} C. Pugh, {\it The $C^{1+\alpha}$ hypothesis in Pesin theory,
} Publ. Math., Inst.
Hautes tud. Sci. 59 (1984), 143-161.


\bibitem{RW} D. Richeson and J. Wiseman, {\it Chain recurrence rates and topological entropy,} Topology
Appl. 156 (2008), 251-261.




\bibitem{Ruelle01} D. Ruelle, {\it Historic behaviour in smooth dynamical systems,} Global Analysis of Dynamical
 Systems (H. W. Broer, B. Krauskopf, and G. Vegter, eds.), Bristol: Institute of Physics
 Publishing, 2001.

\bibitem{Sak} K. Sakai, {\it$C^1$
-stably shadowable chain components,} Ergodic Theory Dyn. Syst. 28 (2008),
987-1029.

\bibitem{Sig} K. Sigmund, {\it Generic properties of invariant measures for axiom A diffeomorphisms,} Invention Math. 11 (1970), 99-109.




\bibitem{Takens} F. Takens,  {\it Orbits with historic behaviour, or non-existence of averages,}  Nonlinearity,  21
                   (2008), 33-36.


\bibitem{TV03}
Takens F, Verbitskiy E. {\it On the variational principle for the topological entropy of certain non-compact sets}. Ergodic Theory $\&$ Dynamical Systems, 2003, 23(1):317-348.

\bibitem{To2010} D. Thompson, {\it The irregular set for maps with the specification property has full topological pressure,}
Dyn. Syst. 25 (2010), no. 1, 25-51.

\bibitem{Tho2012} D. Thompson,  {\it Irregular sets, the $\beta$-transformation and the almost specification property}, Transactions of the American Mathematical Society, 2012, 364 (10): 5395-5414.

{ \bibitem{TV} X. Tian, P. Varandas, {\it Topological entropy of level sets of empirical measures for non-uniformly
expanding maps,}  Discrete and Continuous Dynamical Systems - Series A, 37:10 (2017) 5407-5431.}

\bibitem{Walter} P. Walters,  {\it An introduction to ergodic theory,}
Springer-Verlag, 2001.

\bibitem{WGW} X. Wen, S. Gan and L. Wen, {\it$C^1$
-stably shadowable chain components are hyperbolic,} J.
Differ. Equations 246 (2009), 340-357.

\bibitem{Yo} K. Yano. {\it A remark on the topological entropy of homeomorphisms.} Invent. Math. 59 (1980), 215-220.

\end{thebibliography}
\end{document}